**On the Coefficients of Primitive Normal Polynomials**
N. A. Carella, October 2005.

***Abstract***: The previous paper of the Winter 1998/1999 proved the existence of primitive polynomials, and primitive normal polynomials of degree $n$ with $k$ prescribed coefficients in the finite field $\mathbf{F}_q$ for all sufficiently large powers $q = p^v$ such that $k < p$, and $k < n(1 - \varepsilon) - 1$, $\varepsilon > 0$. This paper presents a longer version of the result on primitive normal polynomials.

**1 Introduction**

Let $q = p^v$ be a sufficiently large prime power, and let $k$ and $n$ be a pair of integers such that $k < p$, and $k < n(1 - \varepsilon) - 1$, $\varepsilon > 0$. This work is a longer version of the analysis on the distribution of the coefficients of primitive polynomials and primitive normal polynomials in finite fields started in [4]. The previous results proved the existence of primitive polynomials and to primitive normal polynomials $f(x) = x^n + a_1 x^{n-1} + \cdots + a_{n-1} x + a_n \in \mathbf{F}_q[x]$ with $k$ prescribed consecutive coefficients $a_1 \neq 0$, $a_2$, ..., $a_k \in \mathbf{F}_q$ in finite fields of odd characteristic $p$. A *primitive polynomial* has roots of multiplicative order $q^n - 1$; and a *normal polynomial* has roots of additive order $x^n - 1$, which is the same as having nonzero trace and linearly independent roots.

***Theorem 1.*** (Extended Coefficient Theorem) Let $q = p^v$, $p$ prime, and let $k$ and $n$ be a pair of integers such that $k < p$, and $k < n(1 - \varepsilon) - 1$, $\varepsilon > 0$. Then there exists a primitive normal polynomial $f(x) \in \mathbf{F}_q[x]$ with $k$ prescribed coefficients for all odd $q \geq q_0$.

The case $k = 1$ was considered in [3], it showed that the existence of primitive normal polynomials of prescribed trace $a_1 \neq 0$ for all pairs $(n, q)$ but a finite numbers of exceptions.

A *primitive normal basis* is a basis $\{\eta, \eta^q, \eta^{q^2}, ..., \eta^{q^{n-1}}\}$ of the vector space $\mathbf{F}_{q^n}$ over $\mathbf{F}_q$ generated by a primitive normal element $\eta \in \mathbf{F}_{q^n}$. The asymptotic proof of the *Primitive Normal Basis Theorem* was first established by both [5], and [8]. And the complete version for all pair $n, q$ was established by [15].

**Theorem 2.** (*Primitive Normal Basis Theorem*)  Let $\mathbf{F}_{q^n}$ be an $n$-degree extension of $\mathbf{F}_q$. Then $\mathbf{F}_{q^n}$ has a primitive normal basis over $\mathbf{F}_q$.

The next result is a refinement of the Theorem 2, it calls for primitive normal elements of arbitrary traces.

**Theorem 3.** (*Primitive Normal Basis Theorem Of Arbitrary Trace*)  For every $a \neq 0$ in $\mathbf{F}_q$, there exists a primitive normal element in $\mathbf{F}_{q^n}$ of trace $a$.

The refinement was proposed as a conjecture in [17], and its first asymptotic proof was completed in [3]. About two years later it was extended to all pairs $n$, $q$ in [7]. A much more general extension remains an open problem.

**Conjecture 4.** (*Morgan-Mullen* 1996)  Let $q$ be a prime power and let $n > 3$ be an integer. Then there exists a completely normal primitive polynomial of degree $n$ over $\mathbf{F}_q$.

## 2 Auxiliary Results

This section introduces several concepts used to study the distribution of elements and the coefficients of polynomials in finite fields.

*Characteristic Functions*

A characteristic function encapsulates certain properties of a subset of elements of $\mathbf{F}_{q^n}$. It effectively filters out those elements that do not satisfy the constraints. The equation of a characteristic function is of the form

$$C(\alpha) = \begin{cases} 1 & \text{if } \alpha \text{ satisfies the properties,} \\ 0 & \text{otherwise,} \end{cases} \tag{1}$$

for all elements $\alpha \in \mathbf{F}_{q^n}$.

*The Characteristic Function of Primitive Elements*

Let $\chi$ be a multiplicative character of order $d = ord(\chi)$, $d \mid q^n - 1$, on $\mathbf{F}_{q^n}$, see [16]. The characteristic function of primitive elements in $\mathbf{F}_{q^n}$ is defined by

$$C_P(\alpha) = \frac{\varphi(q^n - 1)}{q^n - 1} \sum_{d \mid q^n - 1} \frac{\mu(d)}{\varphi(d)} \sum_{ord(\chi) = d} \chi(\alpha), \tag{2}$$

where $\alpha \in \mathbf{F}_{q^n}$, and the arithmetic functions $\mu$ and $\varphi$ are the Mobius and Euler functions on the ring of integers $\mathbb{Z}$ respectively, see [1], and [18]. A product version of this formula has the shape

$$C_P(\alpha) = \frac{\varphi(q^n - 1)}{q^n - 1} \prod_{p \mid q^n - 1} \left(1 - \frac{1}{p-1} \sum_{ord(\chi)=p} \chi(\alpha)\right), \qquad (3)$$

where $p$ runs through the prime divisors of $q^n - 1$. The transformation used to obtain this form is straightforward and appears throughout the literature, other forms of this characteristic function are also possible.

The function $C_P(\alpha)$ is one of the basic tools used in the investigation of the distribution of primitive elements in finite fields. Typical applications are illustrated in [14], [20], [21], etc.

*The Characteristic Function of Normal Elements*

Let $\psi$ be an additive character of order $f(x) = Ord(\psi) \mid x^n - 1$ on $\mathbf{F}_{q^n}$, see [5], and [15]. The characteristic function of normal elements in $\mathbf{F}_{q^n}$ is defined by

$$C_N(\alpha) = \frac{\Phi(x^n - 1)}{q^n} \sum_{d(x) \mid x^n - 1} \frac{M(d(x))}{\Phi(d(x))} \sum_{Ord(\psi)=d(x)} \psi([(x^n - 1)/d(x)] \circ \alpha), \qquad (4)$$

where $\alpha \in \mathbf{F}_{q^n}$, and the arithmetic functions M and $\Omega$ are the Mobius and Euler functions on the ring of polynomials $\mathbf{F}_q[x]$ respectively, see [9], [15], etc, for more details. The product version of this formula has the form

$$C_N(\alpha) = \frac{\Phi(x^n - 1)}{q^n} \prod_{f(x) \mid x^n - 1} \left(1 - \frac{1}{q^{\deg(f(x))} - 1} \sum_{Ord(\psi)=f(x)} \psi([x^n - 1)/f(x)] \circ \alpha)\right), \qquad (5)$$

where $f(x)$ runs through the irreducible factors of $x^n - 1$, see [15], etc.

**Example 5.** For the parameter $n = p^u$, $q = p^v$, the polynomial $x^n - 1 = (x - 1)^n \in \mathbf{F}_q[x]$ and the expression $[(x^n - 1)/g(x)] \circ \alpha = Tr(\alpha)$ is the trace $Tr : \mathbf{F}_{q^n} \rightarrow \mathbf{F}_q$, (since $f(x) = x - 1$ is the only irreducible factor of $x^n - 1$). Thus the characteristic function of normal elements in $\mathbf{F}_{q^n}$ is

$$C_N(\alpha) = (1 - 1/q)\left(1 - \frac{1}{q-1}\sum_{Ord(\psi)=x-1}\psi(Tr(\alpha))\right) = \begin{cases} 0 & \text{if } Tr(\alpha) = 0, \\ 1 & \text{if } Tr(\alpha) \neq 0. \end{cases} \qquad (6)$$

### The Characteristic Function of Primitive Normal Elements

The product of the characteristic function of primitive elements (2) and the characteristic function of normal elements (4) in $\mathbf{F}_{q^n}$ realizes the characteristic function of primitive normal elements. Specifically

$$C_{PN}(\alpha) = \frac{\varphi(q^n-1)}{q^n-1}\frac{\Phi(x^n-1)}{q^n}\sum_{d|q^n-1}\frac{\mu(d)}{\varphi(d)}\sum_{f(x)|x^n-1}\frac{M(f(x))}{\Phi(f(x))}\sum_{Ord(\psi)=f(x)}\sum_{ord(\chi)=d}\chi(\alpha)\psi(\beta), \quad (7)$$

where $\alpha \in \mathbf{F}_{q^n}$, and $\beta = [(x^n-1)/f(x)] \circ \alpha$.

The function $C_{PN}(\alpha)$ is one of the basic tools used in the investigation of the distribution of primitive normal elements in finite fields. Typical applications are illustrated in [17], [4], etc, and in this paper.

### The Characteristic Function of Completely Normal Elements

An element $\eta \in \mathbf{F}_{q^n}$ is *completely normal* if and only if $\eta$ is a normal element in $\mathbf{F}_{q^n}$ over $\mathbf{F}_{q^d}$ for all $d \mid n$. The characteristic function of completely normal elements in $\mathbf{F}_{q^n}$ is constructed from the product of the characteristic functions of normal elements in $\mathbf{F}_{q^n}$ over $\mathbf{F}_{q^d}$ for all $d \mid n$. Here $\mathbf{F}_{q^n}$ is an extension of $\mathbf{F}_{q^d}$ of degree $e = [\mathbf{F}_{q^n} : \mathbf{F}_{q^d}]$, with $n = de$. Thus it follows that the characteristic function of completely normal elements is the product of the individual functions:

$$\begin{aligned} C_{CN}(\alpha) &= \prod_{e|n}\left(\frac{\Phi(x^e-1)}{q^e}\sum_{f(x)|x^e-1}\frac{M(f(x))}{\Phi(f(x))}\sum_{Ord(\psi)=f(x)}\psi([x^e-1]/f(x)]\circ\alpha)\right) \\ &= \prod_{e|n}\prod_{f(x)|x^e-1}\frac{\Phi(x^e-1)}{q^e}\left(1 - \frac{1}{q^{\deg(f(x))}-1}\sum_{Ord(\psi)=f(x)}\psi([x^e-1]/f(x)]\circ\alpha)\right), \end{aligned} \qquad (8)$$

where $\alpha \in \mathbf{F}_{q^n}$, and $f(x)$ runs through the irreducible factors of $x^e - 1 \in \mathbf{F}_{q^d}[x]$.

**Example 6.** For the parameter $n = p^u$, $q = p^v$, the polynomial $x^n - 1 = (x-1)^n \in \mathbf{F}_{q^{p^i}}[x]$ and the expression $[(x^n-1)/f(x)]\circ\alpha = Tr_i(\alpha)$ is the trace $Tr_i : \mathbf{F}_{q^{p^u}} \rightarrow \mathbf{F}_{q^{p^i}}$, $0 \le i < u$, (since $f(x) = x - 1$ is the only irreducible factor of $x^n - 1$). Thus the characteristic function of completely normal elements in $\mathbf{F}_{q^n}$ is

$$C_{CN}(\alpha) = \prod_{i=0}^{u-1}\left(1-1/q^{p^i}\right)\left(1-\frac{1}{q^{p^i}-1}\sum_{Ord(\psi)=x-1}\psi(Tr(\alpha))\right). \tag{9}$$

*Some Probabilities Formulae*

The probabilities $P_1 = \mathrm{P}(\mathrm{ord}(\alpha)=q^n-1)$ and $P_2 = \mathrm{P}(\mathrm{Ord}(\alpha)=x^n-1)$ respectively of primitive elements and normal elements $\alpha$ in a finite field extension $\mathbf{F}_{q^n}$ of $\mathbf{F}_q$ are given by

$$P_1 = \frac{\varphi(q^n-1)}{q^n-1} = \prod_{p\,|\,q^n-1}(1-1/p) \geq \frac{e^{-\gamma}}{\log(q^n-1)} \tag{10}$$

where $p$ ranges over the prime divisors of $q^n - 1$, and

$$P_2 = \frac{\Phi(x^n-1)}{q^n} = \prod_{f(x)\,|\,x^n-1}(1-1/p^{\deg(f)}) \geq (1-1/q)^n, \tag{11}$$

where $f(x)$ ranges over the irreducible divisors of $x^n - 1$ respectively. The degree $\deg(f) = d$ of each irreducible factor $f(x)$ is a divisor of $n$.

The distribution of primitive normal elements is more intricate than either the distribution of primitive elements or the distribution of normal elements. An exact closed form formula for the number of primitive normal bases of $\mathbf{F}_{q^n}$ over $\mathbf{F}_q$ appears to be unknown. However there are asymptotic approximations. For example, the probabilities $P_3 = \mathrm{P}(\mathrm{ord}(\alpha)=q^n-1$ and $\mathrm{Ord}(\alpha)=x^n-1)$ of primitive normal elements is approximated by

$$P_3 \approx \frac{\varphi(q^n-1)}{q^n-1}\frac{\Phi(x^n-1)}{q^n} \geq \frac{e^{-\gamma}(1-1/q)^n}{\log(q^n-1)}. \tag{12}$$

It is clear that if $q$ is sufficiently large, $P_1$ and $P_3$ are essentially the same. In fact the probabilities $P_1$ and $P_2$ are asymptotically independent. This means that the cardinalities of the sets of primitive polynomials and primitive normal polynomials, namely,

$$\frac{\varphi(q^n-1)}{n} \sim \frac{e^{-\gamma}(q^n-1)}{n\log(q^n-1)} \quad \text{and} \quad \frac{\varphi(q^n-1)\Phi(x^n-1)}{n} \sim \frac{e^{-\gamma}(q-1)^n}{n\log(q^n-1)}, \tag{13}$$

are very close. Thus the results for the coefficients of primitive polynomials and primitive normal polynomials are about the same. Another approximation due to [5] for the number of primitive normal elements is

$$PN_n(q) = \frac{\varphi(q^n - 1)\Phi(x^n - 1)}{q^n} + O(q^{(.5+\varepsilon)n}), \tag{14}$$

for all $\varepsilon > 0$.

*Estimate of exponential sums*
A few estimates are required in order to derive nontrivial results on the distribution of elements in finite fields.

**Theorem 7.** If $\chi \neq 1$ is a nontrivial multiplicative character on $\mathbf{F}_{q^n}$, then

$$\left| \sum_{Tr(\xi) \neq 0} \chi(\xi) \right| \leq q^{(n-1)/2}. \tag{15}$$

**Theorem 8.** Let $\psi$ and $\chi$ be a pair of nontrivial additive and multiplicative characters on $\mathbf{F}_{q^n}$, and let $f(x), g(x) \in \mathbf{F}_q[x]$ be polynomials of degrees $k$ and $m$ respectively. Then

$$\left| \sum_{x \in \mathbf{F}_{q^n}} \psi(f(x))\chi(g(x)) \right| \leq (k + m - 1)q^{n/2}, \tag{16}$$

where $f(x), g(x)$ are not $q$ powers and $k + m$ is the number of distinct roots of $f$ and $g$ in the splitting field, see [22], [24], and [25].

*The point counting function*
The point counting function

$$S_a(\xi) = \frac{1}{q} \sum_{\xi \in \mathbf{F}_{q^n}} \sum_{x \in \mathbf{F}_q} \psi(x(f(\xi) - a)) \tag{17}$$

enumerates the cardinality of the solution set $\{ \xi \in \mathbf{F}_{q^n} : Tr(f(\xi)) - a = 0 \}$ of the equation $Tr(f(\xi)) - a = 0$ in $\mathbf{F}_{q^n}$, where $f(x)$ is a function on $\mathbf{F}_{q^n}$, and $a \in \mathbf{F}_q$ is a constant.

*Newton identities in finite fields*
The coefficients of the polynomial $f(x) = x^n + a_1 x^{n-1} + \cdots + a_{n-1}x + a_n$ are given by

$$a_1 = -\sum_{1 \leq i \leq n} z_i, \ \ a_2 = \sum_{1 \leq i_1 < i_2 \leq n} z_{i_1} z_{i_2}, \ \ ..., \ \ a_i = (-1)^i \sum_{1 \leq i_1 < \cdots < i_i \leq n} z_{i_1} \cdots z_{i_i}, \ \ ..., \ \ a_n = (-1)^n \prod_{1 \leq i \leq n} z_i, \tag{18}$$

where $z_1, z_2, \ldots, z_n$ are its roots. The symmetric functions $\sigma_i(z_1,\ldots,z_n)$ and the coefficients $a_i$ are equal up to a sign, that is, $a_i = (-1)^i\sigma_i(z_1,\ldots,z_n)$. The associated power sums are defined by

$$w_i = \sum_{1 \le j \le n} z_j^i = z_1^i + z_2^i + \cdots + z_n^i, \quad i \ge 1. \tag{19}$$

The cyclic structure of the roots of polynomials in cyclic extensions can be utilized to transform the power sums into different forms. Specifically, in finite fields the power sums become

$$w_i = \sum_{1 \le j \le n} z_j^i = \sum_{1 \le j \le n} \alpha^{iq^j} = Tr(\alpha^i), \quad i \ge 1, \tag{20}$$

where $z_i = \alpha^{q^j}$ some $j = 0, 1, \ldots, n - 1$, and $f(\alpha) = 0$. Replacing the power sums in Newton identities (1707) yields the corresponding identities for finite fields:

(1) $\quad Tr(\alpha^k) = -ka_k - \sum_{i=1}^{k-1} a_i Tr(\alpha^{k-i}), \quad 1 \le k \le n,$ (21)

(2) $\quad Tr(\alpha^k) = -\sum_{i=1}^{k-1} a_i Tr(\alpha^{k-i}), \quad n < k.$

These formulae are employed to solve problems about the coefficients of polynomials via the power sums. Fast algorithms for computing the coefficients from the power sums and conversely in $O(n\log(n))$ operations are discussed in [2]. A few of the coefficients are given here in terms of the power sums $w_i = \mathrm{Tr}(\alpha^i)$.

(1) $\quad a_1 = -w_1 = -Tr(\alpha) = -(\alpha + \alpha^q + \alpha^{q^2} + \cdots + \alpha^{q^{n-1}}),$ (22)

(2) $\quad a_2 = (2!)^{-1}(w_1^2 - w_2),$

(3) $\quad a_3 = (3!)^{-1}(-w_1^3 - w_1 w_2 + 2w_3),$

(4) $\quad a_4 = (4!)^{-1}(2w_1^4 + 4w_1^2 w_2 - 3w_2^2 + 4w_1 w_3 - 6w_4),$

(k) $\quad a_k = (k!)^{-1} \sum_e a_e w_1^{e_1} w_2^{e_2} \cdots w_k^{e_k}, \quad e = (e_1,\ldots,e_k), e_i \ge 0, \ k \le n.$

These formula are defined in characteristic $p > k$. A direct calculation of the second coefficient $a_2 = \sum_{0 \le i < j < n} \alpha^{q^i + q^j} = (2!)^{-1}(Tr(\alpha)^2 - Tr(\alpha^2))$ in characteristic $p > 2$ can be accomplished in a few lines. But direct calculations of the other coefficients appear to be a lengthy and difficult task. Even the third coefficient $a_3 = \sum_{0 \le i < j < k < n} \alpha^{q^i + q^j + q^k} = (3!)^{-1}(2Tr(\alpha)^3 - Tr(\alpha)Tr(\alpha^2) - Tr(\alpha^3))$ in characteristic $p > 3$ is difficult to compute directly. But using Newton identities in finite fields (21) this task is a simple algebraic manipulation.

*Coefficients system of equations*

Let $c_1, c_2, ..., c_k \in \mathbf{F}_q$ be constants. The strategy of the analysis is to show that the system of equations

$$Tr(x) = c_1, Tr(x^2) = c_2, ..., Tr(x^k) = c_k \tag{23}$$

has at least one primitive normal element solution $x = \alpha \in \mathbf{F}_{q^n}$. This in turn implies the existence of at least one primitive normal polynomial $f(x)$ with $k$ prescribed consecutive coefficients. The constants $c_1, c_2, ..., c_k$ indirectly prescribes the $k$ coefficients, namely,

$$a_1 = -Tr(\alpha), a_2 = a_2(Tr(\alpha), Tr(\alpha^2)), ..., a_k = a_k(Tr(\alpha), Tr(\alpha^2), ..., Tr(\alpha^k)). \tag{24}$$

The derivation of (24) from (23), which is a direct consequence of Newton identities in finite fields. The application of (23), (24), and its use in the exponential sum (17) are the key ideas in the proofs of several problems on the existence and distribution of the coefficients of certain polynomials over finite fields. I observed the relationship between the trace function and Newton identities in the early 1990's while investigating the trace representations of sequences, and applied it to the theory of coefficients of primitive polynomials in [4].

## 3 The Extended Coefficient Theorem

There are various ways of extending the result on primitive polynomials of degree $n$ with $k$ prescribed coefficients in the finite field $\mathbf{F}_q$ given in [4]. The refinement given here extends it to primitive normal polynomials. The constraint of linearly independent roots on the roots of polynomials restricts the set of primitive normal polynomials to a proper subset of set primitive polynomials. There is one exception; the set of quadratic primitive polynomials and the set of quadratic primitive normal polynomials coincide.

The first line of the proof given below, using a reductio ad absurdum argument, combines the various characteristic functions and the point counting function:

$$\sum_{Tr(\alpha) \neq 0, \ \alpha \in \mathbf{F}_{q^n}} C_N(\alpha) C_P(\alpha) \prod_{i=1}^{k} S_{c_i}(\alpha) = 0 . \tag{25}$$

This claims that there is no primitive normal element solution of equation (23). The usual and standard strategy of dealing with this type of equation is to decompose it into several sums, and then compute a lower estimate, which contradicts the equality, consult [3], [4], [5], [8], [13], [14], [15], [16], [20], and [21] for background details and other references on this type of analysis.

*Proof of Theorem* 1: The number of primitive normal polynomials with $k$ prescribed consecutive coefficients $a_1 \neq 0$, $a_2$, ..., $a_k \in \mathbf{F}_q$, $q$ odd, is given by

(26)

$$N^*(n, q, c_1, ..., c_k) = \sum_{Tr(\xi) \neq 0, \xi \in \mathbf{F}_{q^n}} \left( \frac{\varphi(q^n - 1)}{q^n - 1} \sum_{d | q^n - 1} \frac{\mu(d)}{\varphi(d)} \sum_{ord(\chi) = d} \chi(\xi) \right)$$

$$\times \left( \frac{\Phi(x^n - 1)}{q^n - 1} \sum_{f | x^n - 1} \frac{M(f)}{\Phi(f)} \sum_{Ord(\psi) = f} \psi(Tr(\xi)) \right) \left( \prod_{i=1}^{k} \frac{1}{q} \sum_{x_i \in \mathbf{F}_q} \psi(x_i(Tr(\xi) - c_i)) \right),$$

One possible approach to deal with (26) is to decompose it according to the order $\text{Ord}(\psi) = 1$ or $\text{Ord}(\psi) \neq 1$ of the additive character $\psi$. To accomplish this, rewrite it as

(27)

$$N^*(n, q, c_1, ..., c_k) = \frac{P_1 P_2}{q^k} \sum_{f | x^n - 1} \frac{M(f)}{\Phi(f)} \sum_{Ord(\psi) = f} \sum_{d | q^n - 1} \frac{\mu(d)}{\varphi(d)}$$

$$\times \sum_{ord(\chi) = d} \sum_{x_i \in \mathbf{F}_q} \sum_{Tr(\xi) \neq 0, \xi \in \mathbf{F}_{q^n}} \chi(\xi) \psi(Tr(\xi + \sum_{i=1}^{k} x_i \xi^i) - \sum_{i=1}^{k} c_i x_i).$$

Now separate (27) into two terms corresponding to $\text{Ord}(\psi) = 1$ or $\text{Ord}(\psi) \neq 1$, and simplify to obtain

(28)

$$N^*(n, q, c_1, ..., c_k) = P_1 P_2 \left( q^{n-1}(q - 1) - \sum_{Tr(\xi) \neq 0, \xi \in \mathbf{F}_{q^n}} \chi(\xi) \right) - P_2 N(n, q, c_1, ..., c_k)$$

$$= P_1 P_2 \left( q^{n-1}(q - 1) - \frac{1}{P_1} N(n, q, c_1, ..., c_k) - \sum_{Tr(\xi) \neq 0, \xi \in \mathbf{F}_{q^n}} \chi(\xi) \right),$$

where $q^{n-1}(q-1)$ is the number of elements in $\mathbf{F}_{q^n}$ of nonzero traces, and $N(n, q, c_1, ..., c_k)$ denotes the total number of primitive polynomials $f(x) = x^n + a_1 x^{n-1} + \cdots + a_{n-1} x + a_n \in \mathbf{F}_q[x]$ with $k$ prescribed consecutive coefficients $a_1 \neq 0$, $a_2$, ..., $a_k \in \mathbf{F}_q$.

Replacing the estimates $q^k \leq N(n, q, c_1, ..., c_k) \leq q^{n-s}$, $\prod_{p | q^n - 1} (1 - 1/p)^{-1} = P_1^{-1} < c_\varepsilon q^{n\varepsilon}$, where $s > 0$ is a small integer, $\varepsilon > 0$, and $c_\varepsilon$ is a constant, and the exponential sum estimate yields

$$N^*(n, q, c_1, ..., c_k) \geq P_1 P_2 q^{n-1} \left( q - 1 - \frac{1}{q^{(n-1)/2}} - \frac{c_\varepsilon}{q^{s - \varepsilon n - 1}} \right). \qquad (29)$$

Since the product of the probabilities is in the range $0 < P_1 P_2 < 1$, it readily follows that if $k < p$, then there exists a constant $q_1$ such that $N^*(n,q,c_1,\ldots, c_k) \geq 1$ for all $q \geq q_1$, and $s < \varepsilon n + 1$. ∎

This approach leads to a relatively easy proof of the Theorem 1, it avoids the need to deal with the estimates of various exponential sums, and the function $\Omega(x^n-1)$ which enumerates the distinct irreducible factors of $x^n - 1 \in \mathbf{F}_q[x]$. The case $k = 1$ reduces to the primitive normal basis theorem of arbitrary trace. In fact this is a simpler proof than the first proof given in [3] using a variation of this technique.

## 4 Update

All the published proofs on the existence of primitive polynomials and primitive normal polynomials with $k$ prescribed coefficients assume that the inverse of the term $k!$ must be defined in the finite field. This restricts the value of $k$ to a certain range. The latest papers (see the literature) employ elaborate analysis to overcome this obstacle. Despite these efforts, there are still some values of $k$ that cannot be handle using those methods.

However, the proofs of the existence of primitive polynomials and primitive normal polynomials with $k$ prescribed coefficients are *existence results* and do not require the actual computations of the polynomials. In particular, it is irrelevant whether or not the inverse of the term $k!$ is defined in the finite field. The existence of a solution $x = \alpha \in \mathbf{F}_{q^n}$ of the system of equations

$$Tr(x) = c_1, \, Tr(x^2) = c_2, \, \ldots, \, Tr(x^k) = c_k, \tag{30}$$

where $c_i \in \mathbf{F}_q$, is sufficient. This key observation proves the following.

***Theorem 9.*** (Coefficient Theorem) Let $q = p^v$, $p$ prime, and let $k$ and $n$ be a pair of integers such that $k < n - s$, where $s < n/2$, $s > 0$ is a small fixed integer. Then there exists a primitive normal polynomial $f(x) \in \mathbf{F}_q[x]$ with $k$ prescribed coefficients for all odd $q \geq q_0$.

The small fixed integer $s > 0$ is probably less than 5, and could be a universal constant for all finite fields of sufficiently large cardinality $q$.

Accordingly, for all sufficiently large $q$, almost every coefficient of primitive polynomials and primitive normal polynomials $f(x) \in \mathbf{F}_q[x]$ but a few can be prescribed in advanced. For example, there are primitive polynomials with any consecutive coefficients $a_i$, $a_{i+1}$ prescribed in advanced. This might also be possible for primitive normal polynomials.

**Note:** The proof of Theorem 1 given here is the original one as presented at a conference in 1999. The complete paper was submitted to the journal Finite Fields and Applications around 1998/1999. The results in the paper predate the new crop of papers on the coefficients of primitive polynomials that have been published in the last few years.

**REFERENCES:**

[1] T. M. Apostol, **Introduction to Analytic Number Theory**, Springer-Verlag, N.Y., 1984.

[2] D. Bini, V. Y. Pan, **Polynomial and matrix computations. Vol. 1. Fundamental algorithm**s. Birkhauese Inc., Boston, MA, 1994.

[3] N. A. Carella, *On primitive normal elements of arbitrary traces*, Preprint, Submitted to Finite Fields and Applications in 1996.

[4] __________, *On the coefficients of primitive polynomials*, Preprint, Submitted to Finite Fields and Applications in 1998/1999.

[5] L. Carlitz, *Distribution of Primitive Roots in Finite Fields*, Quarterly J. Math. 4 (1953) p.125-156.

[6] S.D. Cohen, *Primitive elements and polynomials with arbitrary trace*, Discrete Math. 11 (1990), 1-7.

[7] S.D. Cohen, D. Hachenberger, *Primitive normal bases with prescribed trace*. Appl. Algebra Engrg. Comm. Comput. 9 (1999), No. 5, 383-403.

[8] H. Davenport, *Bases for finite fields*, J. London Math. Soc. 43, 1968, p.21-39; Vol. 44, 1969, p. 378.

[9] J. von zur Gathen, M. Giesbrecht, *Constructing normal bases in finite fields*, J. Symbolic Computation, (1990) 10, p.547-560.

[10] K. Ham, G. L. Mullen, *Distribution of irreducible polynomials of small degrees over finite fields*, Math. Comp. Vol. 67, No. 221, 1998, pp.337 - 341.

[11] W. B. Han, *The coefficients of primitive polynomials over finite fields*, Math. Comp. Vol. 65, No. 213, 1996, pp. 331 - 340.

[12] Hansen, G. L. Mullen, *Primitive polynomials over finite fields*, Math. Comp. Vol. 59, No. 211, 1992, pp.639 - 643.

[13] T. Helleseth, *On the Covering Radius of Cyclic Codes and Arithmetic Codes*, Disc. Appl. Math. 11 (1985) 157 - 173.

[14] Dieter Jungnickel, Scott A. Vanstone, *On Primitive Polynomials over Finite Fields*, Journal of Algebra 124, 337-353 (1989).

[15] H.W. Lenstra, R.J. Schoof, *Primitive Normal Bases for Finite Fields*; Math. of Computation, Vol. 48, Number 193, January 1987, pp. 217-231.

[16] Rudolf Lidl, Harald Niederreiter, **Finite Fields**, Encyclopedia of Mathematics and its Applications Vol. 20, 1983, Addison-Wesley Publishing Company.

[17] Ilene Morgan, Gary Mullen, *Primitive Normal Polynomials Over Finite Fields*, Math. Comp. Vol. 63, No. 208, October 1994, p.759-765.

[18] Ivan Niven et al., **An Introduction to the Theory of Numbers**, John Wiley and Sons, N.Y., 1991.

[19] O. Moreno, C. Moreno, *Exponential Sums I and Goppa Codes*, Proc. Amer. Math. Soc. Vol. 11 No. 2, Feb. 1991, pp. 523-531.

[20] Oscar Moreno, *On Primitive Elements Of Trace Equal To 1 In GF($2^m$)*, Discrete Math. 41 (1982) p.53-56.

[21] Oscar Moreno, *On the Existence of a Primitive Quadratic of Trace 1 over GF($p^m$)*, J. Combinatorial Theory Series A 51, 104-110 (1989)

[22] F. Pappalardi, Igor E. Shparlinski, *Artin's Conjecture in Functions fields*, Finite Fields And Theirs Applications, 1, 399 - 404, 1995.

[23] Igor E. Shparlinski, **Computations and Algorithmic Problems in Finite Fields**, Kluwer Academic Press, 1992.

[24] Daquing Wan, *Generators and irreducible polynomials over finite fields*, Math. Comp. Vol. 66, No. 219, 1997, pp.1195 - 1212.

[25] A. Weil, *On some exponential sums*, Proc. Nat. Acad. Sciences. 34, pp.204-207.